\documentclass{amsart}

\usepackage{amsfonts}
\usepackage{amssymb}
\usepackage{amsmath}
\usepackage{hyperref}
\usepackage{graphicx}
\usepackage{color}

\newcommand{\half}{\frac{1}{2}}
\newcommand{\thalf}{\tfrac{1}{2}}
\newcommand{\summ}{\mathop{{\sum}^{\star}}}
\newcommand{\sump}{\mathop{{\sum}^{+}}}

\newcommand{\intt}{\int_{-\infty}^{\infty}}

\def\({\left(}
\def\){\right)}

\numberwithin{equation}{section}

\newtheorem{theorem}{Theorem}[section]
\newtheorem*{theorem*}{Remark}

\newtheorem{lemma}[theorem]{Lemma}

\makeatletter
\@namedef{subjclassname@2010}{%
  \textup{2010} Mathematics Subject Classification}
\makeatother

\begin{document}

\title{Nonvanishing of Dirichlet $L$-functions}

\author{Rizwanur Khan}
\address{
Science Program\\ Texas A\&M University at Qatar\\ PO Box 23874\\ Doha, Qatar}
\email{rizwanur.khan@qatar.tamu.edu, trunghieu.ay@gmail.com }

\author{Hieu T. Ngo}


\begin{abstract} 
We show that for at least $3/8$ of the primitive Dirichlet characters $\chi$ of large prime modulus, the central value $L(1/2,\chi)$ does not vanish.
\end{abstract}

\maketitle

%==================================
\section{Introduction}
%==================================

The zeros of $L$-functions on the critical line are as important in number theory as they are mysterious. At the real point on the critical line (the central point), an $L$-function is expected to vanish only for either a good reason or a trivial reason. A good reason is when the central value has some arithmetic significance which explains why it may vanish. For example, the central value of the $L$-function attached to an elliptic curve over a number field is expected to vanish if and only if the elliptic curve has positive rank (according to the Birch and Swinnerton-Dyer conjecture). A trivial reason is when the functional equation implies that the central value is zero. For instance, the $L$-function of any odd Hecke-Maass form has functional equation $L(\half, f)=-L(\half,f)$. In all other cases, the most extensive success in proving the nonvanishing of $L$-functions has been achieved through the use of mollifiers. For notable examples of the mollifier method, see \cite{kowmicvan, micvan2, iwasar2, sou} as well as the works discussed below.

In this paper, we study the classical nonvanishing problem of primitive Dirichlet $L$-functions. It is conjectured that $L(\half,\chi)\neq0$ for every primitive Dirichlet character $\chi$. Consider for each odd prime $p$ the family of $L$-functions
\begin{align*}
\{ L(s, \chi): \chi \text{ is primitive modulo } p \};
\end{align*}
this family has size $p-2$. 
Viewing $L(\half,\chi)$ as a statistical object, we would like to understand its distribution as $p\to \infty$. One way to get a handle on the distribution is through understanding the moments of $L(\half,\chi)$, but currently only moments of small order are known. Nevertheless this is enough to make some progress in the way of proving that a positive proportion of the family is nonvanishing.

Asymptotic expressions for the first and second moments of $L(\half, \chi)$ are well known. By a result of Heath-Brown \cite{Heath-Brown-a}, we have
\begin{align*}
&\frac{1}{p-2}\summ_{\chi \bmod p} L(\thalf, \chi) \sim 1 \\
&\frac{1}{p-2}\summ_{\chi \bmod p} |L(\thalf, \chi)|^2 \sim \log p,
\end{align*}
where $\summ$ restricts the summation to the primitive characters. The discrepancy between the first and second moments indicates fluctuations in the sizes of the central values. Using these moments and the Cauchy-Schwarz inequality, one can only infer that at least $0 \%$ of the family is nonvanishing, since
\begin{align*}
\frac{1}{p-2}\summ_{\substack{\chi \bmod p\\ L(\half, \chi)\neq 0}} 1\ge \frac{| \frac{1}{p-2}\summ_{\chi \bmod p} L(\thalf, \chi) |^2 }{\frac{1}{p-2}\summ_{\chi \bmod p} |L(\thalf, \chi)|^2  } \gg \frac{1}{\log p}.
\end{align*}
The mollifier method is used to remedy this situation. The origin of the method traces back to the works of Bohr and Landau \cite{bohr-landau} and of Selberg \cite{selberg} on zeros of the Riemann zeta function. The starting idea is to introduce a quantity $M(\chi)$, called the ``mollifier'', which, on average, approximates the inverses of the supposedly nonvanishing values $L(\half, \chi)$. The goal is to choose a mollifier such that the mollified first and second moments are comparable; that is,
\begin{align*}
&\frac{1}{p-2}\summ_{\chi \bmod p} L(\thalf, \chi) M(\chi) \asymp 1\\
&\frac{1}{p-2}\summ_{\chi \bmod p} |L(\thalf, \chi)M(\chi) |^2 \asymp 1.
\end{align*}
From this a positive nonvanishing proportion can be inferred:
\begin{align}
\label{mol-ratio} \frac{1}{p-2}\summ_{\substack{\chi \bmod p\\ L(\half, \chi)\neq 0}} 1 \ge \frac{1}{p-2}\summ_{\substack{\chi \bmod p\\ L(\half, \chi)M(\chi) \neq 0}} 1 \ge \frac{| \frac{1}{p-2}\summ_{\chi \bmod p} L(\thalf, \chi) M(\chi) |^2 }{\frac{1}{p-2}\summ_{\chi \bmod p} |L(\thalf, \chi)M(\chi) |^2  } \gg 1.
\end{align}
Balasubramanian and Murty \cite{balmur} were the first to do this; however their mollifier was inefficient and they obtained only a very small positive proportion of nonvanishing. 

Next came the work of Iwaniec and Sarnak \cite{iwasar}, who introduced a systematic technique that has since served as a model for other families of $L$-functions. Iwaniec and Sarnak took the mollifier
\begin{align}\label{eq:IS-mollifier}
M(\chi) = \sum_{m\le M} \frac{y_m \chi(m) }{m^{\half}},
\end{align}
where $M=p^{\theta}$ is the mollifier length and $(y_m)$ is a sequence of real numbers satisfying $y_m\ll p^{\epsilon}$.
They established the asymptotics of the mollified first and second moments for $\theta<\half$ and found that the choice of coefficients which maximizes the ratio in (\ref{mol-ratio}) is essentially
\begin{align}\label{eq:IS-mollifier-coefficient}
y_m = \mu(m) \frac{\log (\frac{M}{m})}{\log M},
\end{align}
yielding a nonvanishing proportion of
\begin{align*}
\frac{1}{p-2}\summ_{\substack{\chi \bmod p\\ L(\half, \chi)\neq 0}} 1 \ge \frac{\theta}{1+\theta}.
\end{align*}
This can be taken as close to $\frac{1}{3}$ as possible on letting $\theta$ approach $\half$. Computing the mollified moments for larger values of $\theta$ would result in a higher proportion of nonvanishing, but this appears to be very difficult to do. The problem seems to have been attempted by Bettin, Chandee, and Radziwi{\l}{\l}. In \cite{betcharad}, these authors solved the parallel problem for the Riemann zeta function, by obtaining the asymptotics as $T\to \infty$ of 
\begin{align*}
\int_T^{2T} |\zeta(\thalf+it)|^2 |\sum_{m\le M} \frac{y_m}{m^{\half+it}}|^2 \ dt,
\end{align*}
where $M=T^\theta$, for values of $\theta$ slightly larger than $\frac{1}{2}$. However with regard to the problem for Dirichlet $L$-functions, the authors remarked, \emph{``Our proof would not extend to give an asymptotic formula in this case, and additional input is needed.''}

Shortly after the work of Iwaniec and Sarnak, in their study of the nonvanishing of high derivatives of Dirichlet $L$-functions, Michel and VanderKam \cite{micvan} used the ``twisted'' mollifier
\begin{align}\label{twisted} 
M(\chi) = \sum_{m\le M} \frac{y_m \chi(m) }{m^{\half}} + \frac{\overline{\tau}_\chi}{p^\half}\sum_{m\le M} \frac{y_m \overline{\chi}(m) }{m^{\half}},
\end{align}
where $M=p^{\theta}$, $y_m$ is as in \eqref{eq:IS-mollifier-coefficient}, and $\tau_\chi$ is the Gauss sum as defined in their paper. Heuristically, this is a better mimic of $L(\half,\chi)^{-1}$ because the approximate functional equation of $L(\half,\chi)$ essentially consists of a sum of two Dirichlet polynomials, one multiplied by a Gauss sum. A similar two-piece mollifier was first used by Soundararajan \cite{sound} in the context of the Riemann zeta function. Michel and VanderKam \cite{micvan} proved for $\theta<\frac{1}{4}$ a nonvanishing proportion of
\begin{align*}
\frac{1}{p-2}\summ_{\substack{\chi \bmod p\\ L(\half, \chi)\neq 0}} 1 \ge \frac{2\theta}{1+2\theta},
\end{align*}
recovering the $\frac{1}{3}$ proportion of Iwaniec and Sarnak \cite{iwasar}. For this method too, computing the mollified moments for larger $\theta$ would result in a higher proportion of nonvanishing. 

The nonvanishing problem was stuck at the proportion $\frac{1}{3}$ for ten years until Bui \cite{bui} dexterously proved a nonvanishing proportion of $0.3411$. His breakthrough was not to increase the length of any existing mollifier but to use an ingenious new two-piece mollifier. Bui \cite[page 1857]{bui} commented that \emph{``There are two different approaches to improve the results in this and other problems involving mollifiers. One can either extend the length of the Dirichlet polynomial or use some ``better'' mollifiers. The former is certainly much more difficult.''} We take the former, more difficult approach.

Our first idea to attack the nonvanishing problem is to increase the length of the Michel-VanderKam mollifier. This may be a somewhat unexpected avenue because previous attempts at lengthening mollifiers has, as far as we are aware, been directed at the Iwaniec-Sarnak mollifier. Our second idea is to establish an estimate for a trilinear sum of Kloosterman sums with general coefficients (Lemma \ref{mainlemma}). To prove this, we appeal to some work of Fouvry, Ganguly, Kowalski and Michel \cite{fgkm}. The authors thereof proved best possible estimates for sums of products of Kloosterman sums to prime moduli by using powerful algebro-geometric methods (this work built on \cite{fmrs} and was later generalized in \cite{foukowmic}). We stress that although the deepest part of our proof comes from \cite{fgkm}, it is not clear how this work is related to the nonvanishing problem. We figure out this relationship.

Before stating our result, it should be said that the works \cite{iwasar, micvan, bui} actually treat general moduli while we are restricting to prime moduli (which is arguably the most interesting case).

\begin{theorem}\label{main} Let $\epsilon>0$ be arbitrary. For all primes $p$ large enough in terms of $\epsilon$, there are at least $\(\frac{3}{8}-\epsilon\)$ of the primitive Dirichlet characters $\chi$ ({\rm mod} $p$) for which $L(\half, \chi)\neq 0$.
\end{theorem}

\noindent The significance of our work is that we show for the first time how to increase the length of a classical mollifier in this context. An interesting open problem that remains is to increase the length of the Iwaniec-Sarnak mollifier. Our nonvanishing proportion $\frac{3}{8}$ improves upon that of Bui for prime moduli. For general moduli, Bui's nonvanishing proportion 0.3411 is still the best known.

Throughout the paper, we use the standard convention that $\epsilon$ denotes an arbitrarily small positive constant which may differ from one occurrence to the next, and that the implied constants in the various estimates depend on $\epsilon$.

%==================================
\section{The work of Michel and VanderKam}
%==================================

We briefly summarize the mollifier method of Michel and VanderKam \cite{micvan}, setting the ground for our further discussion. 

Let the mollifier $M(\chi)$ be given by \eqref{twisted} where the mollifier length is $M=p^\theta$ and the real mollifying coefficients $y_m$ are given by \eqref{eq:IS-mollifier-coefficient}. Michel and VanderKam asymptotically evaluated the mollified first moment
\begin{align*}
\frac{2}{p-2}\sump_{\chi \bmod p} L(\thalf, \chi) M(\chi)
\end{align*}
for $\theta<\half$, where $\sump$ restricts the summation to the even primitive characters. The evaluation for the odd primitive characters is entirely similar. They evaluated the mollified second moment
\begin{align}
\label{2ndmol}  &\frac{2}{p-2}\sump_{\chi \bmod p} |L(\thalf, \chi) M(\chi)|^2\\
\nonumber &= \frac{4}{p-2} \sump_{\chi \bmod p} |L(\thalf, \chi) |^2 \Big|\sum_{m\le M} \frac{y_m \chi(m)}{m^\half} \Big|^2 \\ 
\nonumber &+ \frac{4}{p-2} \sump_{\chi \bmod p} |L(\thalf, \chi) |^2 \frac{\tau_\chi}{p^\half} \Big(\sum_{m\le M} \frac{y_m \chi(m)}{m^\half} \Big)^2
\end{align}
for $\theta<\frac{1}{4}$; see \cite[Equation (10)]{micvan} for the above identity. An asymptotic for the first sum on the right hand side of \eqref{2ndmol} is derived for $\theta<\half$, as was done by Iwaniec and Sarnak \cite{iwasar}, but the second sum is more difficult and could only be handled for $\theta<\frac{1}{4}$. In the end, the main terms of the mollified moments of Michel and VanderKam yield a nonvanishing proportion of $\frac{2\theta}{1+2\theta}$, by taking $P_0(t)=t$ in \cite[section 7]{micvan}.

Let us concentrate on the second sum on the right hand side of (\ref{2ndmol}). Recall the standard approximate functional equation (see for example \cite[Equation (3)]{micvan}):
\begin{align}
\label{afe} |L(\thalf, \chi)|^2= 2\sum_{n_1,n_2\ge 1} \frac{\chi(n_1)\overline{\chi}(n_2)}{(n_1n_2)^\half}V\Big(\frac{n_1n_2}{p}\Big),
\end{align}
where
\begin{align*}
V(x)=\frac{1}{2\pi i} \int_{(2)} \frac{\Gamma(\frac{s}{2}+\frac{1}{4})^2}{\Gamma( \frac{1}{4})^2}(\pi x)^{-s} \frac{ds}{s}.
\end{align*}
By moving the line of integration, one shows that $V(x)\ll_c x^{-c}$ for any $c>0$, whence the sum in \eqref{afe} is essentially supported on $n_1n_2\le p^{1+\epsilon}$. Therefore
\begin{align}
\nonumber&\frac{4}{p-2} \sump_{\chi \bmod p} |L(\thalf, \chi) |^2 \frac{\tau_\chi}{p^\half} \Big(\sum_{m\le M} \frac{y_m \chi(m)}{m^\half} \Big)^2 \\
\label{2nd} &= \sum_{\substack{n_1,n_2\ge 1\\ m_1,m_2\le M}} \frac{y_{m_1}y_{m_2}}{(n_1n_2m_1m_2)^\half} V\Big(\frac{n_1n_2}{p}\Big)\frac{4}{p-2} \sump_{\chi \bmod p} \frac{\tau_\chi}{p^\half} \chi(n_1m_1m_2)\overline{\chi}(n_2).
\end{align}
By \cite[Equation (17)]{micvan} or \cite[Equation (3.4)]{iwasar}, for $(n,p)=1$ we have 
\begin{align*}
\displaystyle \sump_{\chi \bmod p} \tau_\chi \chi(n) = p\cos \Big(\frac{2\pi \, \overline{n}}{p}\Big) + O(1),
\end{align*}
so that (\ref{2nd}) equals
\begin{align}
\label{2moll} \frac{4}{p^\half} {\rm Re} \sum_{\substack{n_1,n_2\ge 1\\ m_1,m_2\le M\\(n_1 n_2 m_1m_2,p)=1}} \frac{y_{m_1}y_{m_2}}{(n_1n_2m_1m_2)^\half} V\Big(\frac{n_1n_2}{p}\Big) e\Big(\frac{ n_2 \, \overline{n_1 m_1 m_2}}{p}\Big) + O\Big(\frac{M}{p^{1-\epsilon}} + p^{-\epsilon}\Big)
\end{align}
for any $\epsilon>0$, where $e(x)=e^{2\pi i x}$ and $\overline{n}$ denotes the multiplicative inverse of $n$ mod $p$ for $(n,p)=1$. The terms with $m_1m_2=1$ contain a main term of (\ref{2nd}); see \cite[section 6]{micvan}. Consider the rest of the terms in dyadic intervals. Let 
\begin{align}  
\label{bdef} &\mathcal{B}(M_1,M_2,N_1,N_2) \\
\nonumber &= \frac{1}{(pM_1M_2N_1N_2)^\half} \sum_{\substack{n_1,n_2\ge 1\\ M_1\le m_1 \le 2M_1\\ M_2 \le m_2 \le 2M_2 \\ (n_1 n_2 m_1m_2,p)=1 }}  y_{m_1}y_{m_2}  e\Big(\frac{ n_2 \overline{n_1 m_1 m_2}}{p}\Big) V\Big(\frac{n_1n_2}{p}\Big) f_1\Big(\frac{n_1}{N_1}\Big) f_2\Big(\frac{n_2}{N_2}\Big)
\end{align}
for $2\le M_1M_2 \le M^2$, $1\le N_1 N_2\le p^{1+\epsilon}$ and any smooth functions $f_1,f_2$ compactly supported on the positive reals. Michel and VanderKam \cite[Equations (24) and (27)]{micvan} proved the bounds
\begin{align}
\label{firstbound} &\mathcal{B}(M_1,M_2,N_1,N_2)\ll p^\epsilon \Big(\frac{M^2N_1}{pN_2}\Big)^\half
\end{align}
and 
\begin{align}
\label{secondbound} &\mathcal{B}(M_1,M_2,N_1,N_2)\ll p^\epsilon \Big(\frac{M^2N_2}{N_1}\Big)^\half + \frac{M}{p^{1-\epsilon}}. 
\end{align}
These bounds together yield $\mathcal{B}(M_1,M_2,N_1,N_2)\ll p^{-\epsilon}$, provided that $M\le p^{\frac{1}{4}-\epsilon}$. Thus the contribution to (\ref{2moll}) of the terms with $m_1m_2\ge 2$ is $O(p^{-\epsilon})$ for $\theta<\frac{1}{4}$.

In the next section we will show how to improve the bound (\ref{secondbound}), in the ranges where (\ref{firstbound}) is not useful. This together with (\ref{firstbound}) will imply that 
$$\mathcal{B}(M_1,M_2,N_1,N_2)\ll p^{-\epsilon}$$ 
for larger values of $\theta$, thereby extending the asymptotics of Michel and VanderKam.

%==================================
\section{Proof of Theorem \ref{main}}
%==================================

To get the bounds (\ref{firstbound}) and (\ref{secondbound}), Michel and VanderKam obtained cancellation in only the $(n_1,n_2)$-sums of $\mathcal{B}(M_1,M_2,N_1,N_2)$. On the other hand, we will use the $(m_1,m_2)$-sums to our advantage. To set up for this, we first prove some estimates for averages of products of Kloosterman sums. Let
\begin{align*}
S(a,b;c)=\sum_{\substack{x \bmod c\\ x\overline{x} \equiv 1 \bmod c}} e\Big(\frac{ax+b\overline{x}}{c}\Big)
\end{align*}
denote the Kloosterman sum.
The following lemma is a consequence of a result of Fouvry, Ganguly, Kowalski and Michel \cite{fgkm}.

\begin{lemma}\label{kloosterman} 
For $B\le p$ we have
\begin{align}
\label{sumprod} \sum_{\substack{b_1,b_2,b_3,b_4\le B\\ (b_1b_2b_3b_4,p)=1}}  \Big| \sum_{h\bmod p} S(h,\overline{b}_1;p) S(h,\overline{b}_2;p) S(h,\overline{b}_3;p) S(h,\overline{b}_4;p) \Big|\ll B^4p^\frac{5}{2} +B^2p^3.
\end{align}
\proof
Write the left hand side of (\ref{sumprod}) as
\begin{align*}
\sum_{\substack{b_1,b_2,b_3,b_4\le B\\ (b_1b_2b_3b_4,p)=1}} = \sum_{\substack{b_1,b_2,b_3,b_4\le B\\ (b_1,b_2,b_3,b_4)\in \mathfrak{D}\\ (b_1b_2b_3b_4,p)=1}} + \sum_{\substack{b_1,b_2,b_3,b_4\le B\\ (b_1,b_2,b_3,b_4)\notin \mathfrak{D}\\ (b_1b_2b_3b_4,p)=1}} 
\end{align*}
where $\mathfrak{D}$
is the set of tuples $(b_1,b_2,b_3,b_4)$ such that no component $b_i$ is distinct mod $p$ from the others. Note that $|\mathfrak{D}|\ll B^2$. 

On the one hand, it follows from the Weil bound for Kloosterman sums that
\begin{align*}
\sum_{\substack{b_1,b_2,b_3,b_4\le B\\ (b_1,b_2,b_3,b_4)\in \mathfrak{D}\\ (b_1b_2b_3b_4,p)=1}}  \Big| \sum_{h\bmod p} S(h,\overline{b}_1;p) S(h,\overline{b}_2;p) S(h,\overline{b}_3;p) S(h,\overline{b}_4;p) \Big|\ll B^2p^3.
\end{align*}
On the other hand,
if $(b_1,b_2,b_3,b_4)\notin \mathfrak{D}$, then in the language of \cite[Definition 3.1]{fgkm}, $(\overline{b}_1,\overline{b}_2,\overline{b}_3,\overline{b}_4)$ is not in ``mirror configuration''. Thus \cite[Proposition 3.2]{fgkm} asserts that
\begin{align*}
 \sum_{h\bmod p} S(h,\overline{b}_1;p) S(h,\overline{b}_2;p) S(h,\overline{b}_3;p) S(h,\overline{b}_4;p) \ll p^\frac{5}{2},
\end{align*}
saving a factor of $p^\half$ over Weil's bound.
So 
\begin{align*}
\sum_{\substack{b_1,b_2,b_3,b_4\le B\\ (b_1,b_2,b_3,b_4)\notin \mathfrak{D}\\ (b_1b_2b_3b_4,p)=1}}  \Big| \sum_{h\bmod p} S(h,\overline{b}_1;p) S(h,\overline{b}_2;p) S(h,\overline{b}_3;p) S(h,\overline{b}_4;p) \Big|\ll B^4p^\frac{5}{2}.
\end{align*}
The lemma follows.
\endproof
\end{lemma}

Let now
\begin{align*}
\mathcal{S}=\sum_{\substack{1\le |n| \le N\\ 1\le a\le A\\ 1\le b\le B}} x_n y_a z_b S(n,\overline{ab};p),
\end{align*}
where the coefficients satisfy $x_n,y_a, z_b \ll p^\epsilon$, $y_a=0$ for $p|a$, and $z_b=0$ for $p|b$. 

\begin{lemma} \label{mainlemma} 
For $NA\le \frac{p}{2}$ and $B\le p$, we have 
\begin{align*}
\mathcal{S} \ll p^{\epsilon} N^\frac{3}{4} A^\frac{3}{4} (B p^\frac{5}{8} + B^\frac{1}{2}p^\frac{3}{4}).
\end{align*}
\end{lemma}
\proof
On applying the Cauchy-Schwarz inequality, we infer
\begin{align*}
\mathcal{|S|}^2\ll p^\epsilon N A \sum_{\substack{|n| \le N\\ a \le A}} 
			\Big| \sum_{b \le B} z_{b} S(n\,\overline{a}, \overline{b};p) \Big|^2.
\end{align*}
Hence 
\begin{align}
\label{nu-h} \mathcal{|S|}^2\ll p^\epsilon NA 
	\sum_{h \bmod p} \nu(h) \Big| \sum_{b\le B} z_b S(h,\overline{b};p) \Big|^2
\end{align}
where
\begin{align*}
\nu(h) = \sum_{\substack{|n| \le N \\ a\le A \\  n \, \overline{a} \equiv h \bmod p}}  1.
\end{align*}

On applying Cauchy-Schwarz to (\ref{nu-h}), we find that 
\begin{align}
\label{s4} \mathcal{|S|}^4\ll p^\epsilon N^2A^2 
	\Big( \sum_{h \bmod p} \nu(h)^2\Big) 
	\Big( \sum_{h \bmod p}  \Big| \sum_{b\le B} z_b S(h,\overline{b};p) \Big|^4 \Big).
\end{align}
Observe that 
\begin{align*}
 \sum_{h \bmod p} \nu(h)^2 = 
	\sum_{\substack{  |n_1|,|n_2|\le N \\ a_1, a_2\le A \\ n_1 \overline{a}_1 \equiv n_2 \overline{a}_2  \bmod p }}  1 =  
	\sum_{\substack{|n_1|,|n_2|\le N \\ a_1, a_2\le A \\  n_1 a_2 \equiv n_2 a_1  \bmod p}} 1.
 \end{align*}
Since $NA\le \frac{p}{2}$ by assumption, it follows that
\begin{align*}
 \sum_{h \bmod p} \nu(h)^2  =  
	\sum_{\substack{  n_1 a_2= n_2 a_1\\ |n_1|,|n_2| \le N\\  a_1, a_2\le A }}  1  \ll p^{\epsilon} N A.
\end{align*}
Therefore (\ref{s4}) becomes
\begin{align*}
\mathcal{|S|}^4\ll p^\epsilon N^3A^3  
	\sum_{\substack{b_1,b_2,b_3,b_4\le B\\ (b_1b_2b_3b_4,p)=1}}  
	\Big| \sum_{h\bmod p} S(h,\overline{b}_1;p) S(h,\overline{b}_2;p) 
		S(h,\overline{b}_3;p) S(h,\overline{b}_4;p) \Big|.
\end{align*}
Finally, we apply Lemma \ref{kloosterman} to conclude that
\begin{align*}
\mathcal{|S|}^4\ll p^\epsilon N^3A^3  (B^4p^\frac{5}{2} +B^2p^3).
\end{align*}
The lemma is proved.
\endproof

We are in a position to prove a new bound for our nonvanishing problem.

\begin{lemma} \label{newbound}
For $\frac{N_1}{N_2} > p^{\epsilon} M$ and $M< p^{\half-\epsilon}$, we have
\begin{align}
\label{thirdbound}  \mathcal{B}(M_1,M_2,N_1,N_2) \ll  p^\epsilon \Big(\frac{N_2M^3}{N_1p^3}\Big)^\frac{1}{4} \Big( p^\frac{5}{8} +\frac{p^\frac{3}{4}}{M^\frac{1}{2}}\Big) +  \frac{M}{p^{\half-\epsilon}}.
\end{align}
\end{lemma}
\proof

In (\ref{bdef}), separate $n_1$ into residue classes modulo $p$ and apply the Poisson summation formula to get
\begin{align}
\nonumber &\mathcal{B}(M_1,M_2,N_1,N_2)  \\
\label{poisson} &= \frac{1}{(pM_1M_2N_1N_2)^\half} \frac{N_1}{p}\sum_{\substack{-\infty<k<\infty\\ n_2\ge 1\\ M_1\le m_1 \le 2M_1\\ M_2 \le m_2 \le 2M_2 \\ (n_2 m_1m_2,p)=1 }}  y_{m_1}y_{m_2} S(k n_2 , \overline{m_1m_2};p) f_2\Big(\frac{n_2}{N_2}\Big) F(k)
\end{align}
where
\begin{align*}
F(k)=  \intt f_1(x) V\Big(\frac{ xN_1n_2}{p}\Big)e\Big(\frac{-xkN_1}{p}\Big) dx.
\end{align*}
Repeatedly integrating by parts, we find that $F(k)\ll_c \big(\frac{kN_1}{p}\big)^{-c}$ for any $c>0$. Thus the $k$-sum may be restricted to $|k|\le \frac{p^{1+\epsilon}}{N_1}$.

The contribution to (\ref{poisson}) of the terms with $k=0$ is 
\begin{align*}
&\frac{1}{(pM_1M_2N_1N_2)^\half} \frac{N_1}{p}\sum_{\substack{ n_2\ge 1\\ M_1\le m_1 \le 2M_1\\ M_2 \le m_2 \le 2M_2 \\ (n_2 m_1m_2,p)=1 }}  y_{m_1}y_{m_2} S(0 , \overline{m_1m_2};p) f_2\Big(\frac{n_2}{N_2}\Big) F(0)\\
&\ll \frac{(N_1N_2M_1M_2)^\half}{p^{1-\epsilon}}\ll \frac{M}{p^{\half-\epsilon}}.
\end{align*}
This is the last term in (\ref{thirdbound}).
The contribution of the terms with $|k|>0$ is bounded using Lemma \ref{mainlemma}, by putting
\begin{align*}
&n=kn_2, \ \ x_n = f_2\big(\tfrac{n_2}{N}\big)F(k) \text{ if } (n_2,p)=1, \ \ x_n = 0 \text{ if } p|n_2, \ \ N=\tfrac{N_2 p^{1+\epsilon}}{N_1} \\
&a=m_1, \ \ y_a= y_{m_1}, \ \ A=2M_1\\
&b=m_2, \ \ z_b = y_{m_2}, \ \ B=2M_2.
\end{align*}
Note that the conditions of Lemma \ref{mainlemma}, namely $B\le p$ and $NA\le \frac{p}{2}$, are satisfied by the assumptions that $M< p^{\half-\epsilon}$ and that $\frac{N_1}{N_2} > p^{\epsilon} M$. The bound (\ref{thirdbound}) follows.
\endproof 

Finally, we sum up the work done to arrive at the following power-saving result.

\begin{lemma}\label{lem:o1}
We have $\mathcal{B}(M_1,M_2,N_1,N_2)\ll p^{-\epsilon}$ for $M<p^{\frac{3}{10}-\epsilon}$.
\proof
Assume first that $M<p^{\frac{1}{3}-\epsilon}$. 
If $\frac{N_1}{N_2} \le p^{\epsilon} M$, it follows from \eqref{firstbound} that $\mathcal{B}(M_1,M_2,N_1,N_2)\ll p^{-\epsilon}$, whence the lemma follows. 

We therefore suppose that $\frac{N_1}{N_2} > p^{\epsilon} M$. Now since the conditions of Lemma \ref{newbound} are met, we have the bound (\ref{thirdbound}). In this bound, we may suppose that $\frac{N_2}{N_1}<\frac{M^2}{p^{1-\epsilon}}$, since otherwise by (\ref{firstbound}), we have $\mathcal{B}(M_1,M_2,N_1,N_2)\ll p^{-\epsilon}$. Thus (\ref{thirdbound}) becomes
\begin{align*}
 \mathcal{B}(M_1,M_2,N_1,N_2) \ll  \frac{M^\frac{5}{4}}{p^{1-\epsilon}} \Big( p^\frac{5}{8} +\frac{p^\frac{3}{4}}{M^\frac{1}{2}}\Big) + p^{-\frac{1}{6}+\epsilon}.
\end{align*}
The bound is $O(p^{-\epsilon})$ precisely when $M\ll p^{\frac{3}{10}-\epsilon}$.
The lemma follows.
\endproof

\proof[Proof of Theorem \ref{main}]
By Lemma \ref{lem:o1}, the nonvanishing proportion $\frac{2\theta}{1+2\theta}$ of Michel and VanderKam is valid for any $\theta < \frac{3}{10}$. On letting $\theta$ approach $\frac{3}{10}$, we infer that the nonvanishing proportion is at least $\frac{3}{8}-\epsilon$ for any $\epsilon>0$.
\endproof

\end{lemma}

\bibliographystyle{amsplain}

\bibliography{dirichletnonvanishing}

\end{document}